# Diffeomorphism Neural Operator for various domains and parameters of partial differential equations

Zhiwei Zhao[†,1], Changqing Liu[†,1], Yingguang Li[*,1], Zhibin Chen[1], Xu Liu[2]

[1]: College of Mechanical and Electrical Engineering, Nanjing University of Aeronautics and Astronautics, Nanjing 210016, China

[2]: School of Mechanical and Power Engineering, Nanjing Tech University, Nanjing 211800, China

[*] Corresponding author.

*E-mail address:* liyingguang@nuaa.edu.cn (Yingguang Li).

[†] These authors contributed equally to this work.

**Abstract:**

In scientific and engineering applications, solving partial differential equations (PDEs) across various parameters and domains normally relies on resource-intensive numerical methods. Neural operators based on deep learning offered a promising alternative to PDEs solving by directly learning physical laws from data. However, the current neural operator methods were limited to solve PDEs on fixed domains. Expanding neural operators to solve PDEs on various domains hold significant promise in medical imaging, engineering design and manufacturing applications, where geometric and parameter changes are essential. This paper presents a novel neural operator learning framework for solving PDEs with various domains and parameters defined for physical systems, named diffeomorphism neural operator (DNO). The main idea is that a neural operator learns in a generic domain which is diffeomorphically mapped from various physics domains expressed by the same PDE. In this way, the challenge of operator learning on various domains is transformed into operator learning on the generic domain. The generalization performance of DNO on different domains can be assessed by a proposed method which evaluates the geometric similarity between a new domain and the domains of training dataset after diffeomorphism. Experiments on Darcy flow, pipe flow, airfoil flow and mechanics were carried out, where harmonic and volume parameterization were used as the diffeomorphism for 2D and 3D domains. The DNO framework demonstrated robust learning capabilities and strong generalization performance across various domains and parameters.

**Main**

Partial differential equations (PDEs) have been used to mathematically represent physical systems in many scientific and engineering applications across different geometric domains and initial/boundary conditions (parameters). Solving PDEs by traditional numerical methods requires considerable computation for recalculations even with minor variations in domains and parameters, which is very time-consuming and resource-intensive [1–3]. Therefore, solving PDEs efficiently and accurately across various domains and parameters has consistently posed significant challenges in scientific research and engineering applications, such as in medical imaging, material science, engineering simulation, where variations in human tissue, material structure, and product shape contribute to not only parameter variations but also domain variations.

Recent advancements in artificial intelligence (AI) including the application of deep learning techniques to efficiently solving complex physical systems described by PDEs, exhibited noteworthy performance in mechanics, thermal transfer, fluid dynamics and climate models[4–11]. A more recently reported approach considered differential operators as mappings between parameter function space and solution function space, and learned the physical laws directly from data in the same class of PDEs with different parameters. The differential operator was named neural operator[12], and reported examples include deep neural operator (DeepONet)[13], Fourier neural operator (FNO)[14], graph kernel network (GKN)[15] and Transformer neural operator[16,17]. In recent years, research efforts had focused on extending the capabilities of neural operators, e.g., exploring complex geometry domain[18,19], multi fields input[20,21], symmetry characteristics[22], physics fields laws[23] and their applications in weather prediction, heterogeneous material modeling and more[24–28]. However, existing neural operator frameworks had a major limitation, i.e., they addressed physical systems with varying parameters but only on fixed domains. Therefore, it always requires re-training even for a minor change in the domain, resulting in high time and training costs. More general frameworks suitable for various domains and parameters had not been thoroughly studied and remains challenging. Although there were attempts to solve PDEs on various domains through transfer learning[29] and expanding to a larger domain[30], there still exist a series of issues to be resolved, such as the requirements of a larger number of labelled data and limited generalization ability.

This paper presents a novel neural operator framework, termed Diffeomorphism Neural Operator (DNO) framework, which was designed to construct a neural operator across various domains and parameters, with these domains being diffeomorphism to a 'generic domain' as shown in (Fig. 1a). The generic domain of a PDE is a new concept proposed in this paper, as opposed to the various physical domains of the same PDE. The underlying principle is that we addressed the challenge of learning function mappings within varied domains by transforming it into the task of learning operators on a generic domain. Diffeomorphism facilitated the mapping of different domains onto a generic domain of the same PDE, along with the parameter function and the solution function defined on them. Then, neural operators learned the mapping between the parameter function and the domain function to the solution function on the generic domain. By learning on the generic domain, the neural operator was able to accurately predict the solution function on various domains. More importantly, diffeomorphism was able to deal with complex domains that were not in the training dataset, thus proving that the trained neural operator had been generalized, i.e., it was able to predict solutions accurately for new domains with different shapes and sizes not even in the training dataset. In summary, the proposed DNO framework extended the scope of previous neural operators which were focused on fixed domains. The framework is also flexible in the selection of diffeomorphism algorithms and neural operator algorithms.

The paper also provided examples on several datasets to demonstrate the advantages of the DNO framework, which proved its remarkable stability in accuracy. As shown in Fig. 1c, for 2D Darcy flow, pipe flow and airfoil flow problems, the harmonic mapping algorithm and FNO were adopted, and for 3D mechanics (stress-deformation) problems, volume parameterization mapping algorithm and FNO were applied. FNO used in the examples was an improved FNO with additional geometric information from domains (Fig. 1b). Furthermore, the examples demonstrated the extensive generalization capability of DNO across domains of varying shapes and sizes. Meanwhile, the generalization performance of DNO was assessed based on a proposed method by evaluating the geometric similarity between any new domains and the training domain on the generic domain.

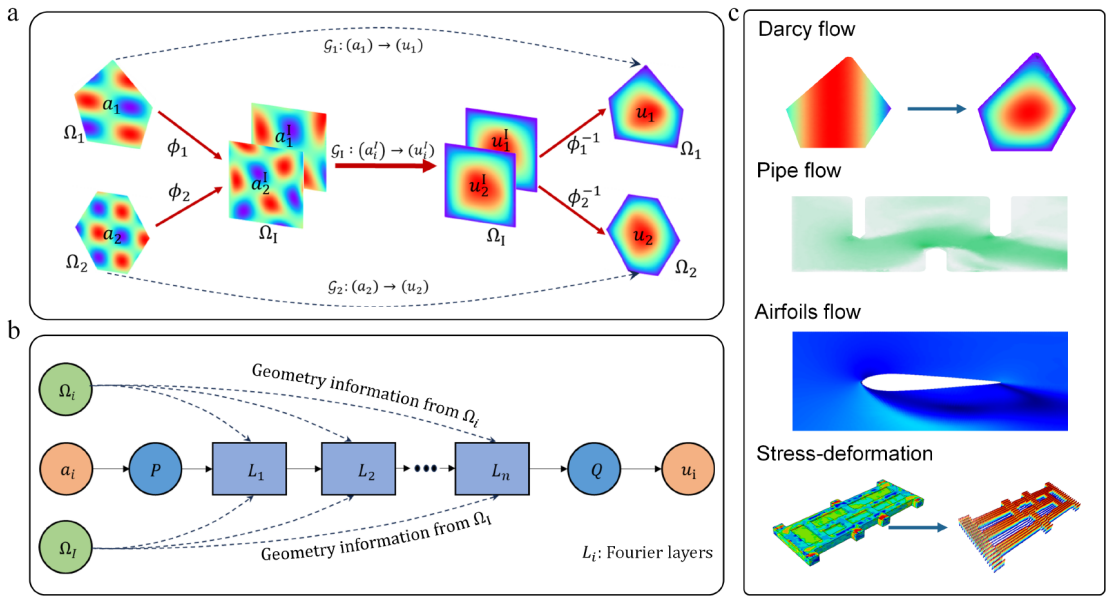

Fig. 1 | The underlying principle of the proposed diffeomorphism neural operator (DNO): (a) The DNO framework. (b) The DNO is an improved FNO with geometry information from both physics domain and generic domain. (c) Examples used to demonstrate advantages of the DNO framework

## Results

### Overview of the proposed diffeomorphism neural operator framework

In a class of PDE problems, the PDE is defined on a set of bounded domains $\mathbf{\Omega} = \{\Omega_1, \dots, \Omega_n\}$ with different shapes. In each domain $\Omega_i \subset \mathbb{R}^d$, there is a varying initial/boundary condition (parameter function) $a_i \in \mathcal{A}$, and corresponding solution function $u_i \in \mathcal{U}$. Assuming that for each domain $\Omega_i$, there is an operator $\mathcal{G}_i: a_i \to u_i$. Here we want to construct an operator $\mathcal{G}_I$ to represent all the individual operators, i.e., $\mathcal{G} = \{\mathcal{G}_1, \dots, \mathcal{G}_n\}$ on domains $\mathbf{\Omega}$, and could also be generalized to new domains which are not in domains $\mathbf{\Omega}$.

The principle of the proposed DNO framework is presented in Fig. 1. As shown in Fig.1a, the proposed framework maps domains $\Omega_i \in \mathbf{\Omega}$, i.e., the corresponding parameter function $a_i$ and solution function $u_i$ of each domain $\Omega_i$, through $\boldsymbol{\phi} = \{\phi_1, \dots, \phi_n\}$ to a generic domain $\Omega_I$ (Fig. 1c), where $\phi_i$ is diffeomorphism from original physics domain $\Omega_i$ to the generic domain $\Omega_I$. Then, the mapping relationships between the parameter functions $\boldsymbol{a}^I = \{a_1^I, \dots, a_n^I\}$ and solution functions $\boldsymbol{u}^I = \{u_1^I, \dots, u_n^I\}$, where $a_i^I = a_i \circ \phi_i$ and $u_i^I = u_i \circ \phi_i$, are learned in the generic domain $\Omega_I$ by neural operator $\mathcal{G}_I$.

$$\mathcal{G}_I : a_i^I \to u_i^I \tag{1}$$

In the proposed DNO framework, the geometric domain information from both the physics domain and the generic domain are integrated, i.e., $\mathcal{G}_I : (a_i^I, G_i, G_I) \to u_i^I$, where $G_i$ and $G_I$ could be the parameterization of the physics and generic domains. Because the diffeomorphism $\boldsymbol{\phi}$ is a continuously differentiable mapping, meaning that both $\phi_i$ and $\phi_i^{-1}$ are continuously differentiable, thus, after training neural operator $\mathcal{G}_I$, $u_i$ could be obtained by:

$$u_i = \mathcal{G}_i(a_i) = u_i^I \circ \phi_i^{-1} = \mathcal{G}_I(a_i^I, G_i, G_I) \circ \phi_i^{-1} \tag{2}$$

Our experimental results shown different impacts on neural operator generalization from different domain information, i.e., geometric information from the generic domain is essential for mesh invariance, while geometric information from the physics domain enhanced the stability of generalization with changes in size and shape (more results of experiments provided in Supplement III). More details about the diffeomorphism algorithm, architecture of neural operators, and the method for computing generic domains are presented in the **Methods** section.

### Darcy flow equations

Darcy's law characterizes the pressure of a fluid flowing through a porous medium with a specified permeability[31]. In this research, a steady state of the Darcy flow equation was considered, and PDE was a second order, linear and elliptic with a Dirichlet boundary. Two dimensional (2D) Darcy flow equations have numerous applications in science and engineering, including modeling the pressure of subsurface flow, the deformation of linearly elastic materials, and the electric potential in conductive materials.

For a domain $\Omega_i \in \mathbf{\Omega}$, the steady state of the 2D Darcy flow equation in $(x, y)$ plane satisfies the following equation.

$$-\nabla\big(a(x,y)\nabla(u(x,y))\big) = F(x,y) \quad (x,y) \in \Omega_i \tag{3}$$
$$u(x,y) = 0 \quad (x,y) \in \partial\Omega_i \tag{4}$$

where $a(x,y)$ is the diffusion coefficient (parameter function) in $(x,y)$ plane, $u(x,y)$ is the solution function and $F(x,y)$ is the force function, which was set as 1 in this case.

Here the goal was to train a neural operator $\mathcal{G}_I$, with input of parameter function $a(x,y)$ and output of solution function $u(x,y)$ on various domains $\mathbf{\Omega}$. The training dataset consisted of pentagon-shaped domains (more dataset details provided in Supplement I), each contained within a $10 \times 10$ area. The corresponding parameter functions were randomly generated according to $a = c_1 \times \sin(x/10) - c_2 \times x \times (x - 10) + 2$, where $c_1, c_2$ followed a uniform distribution between 0.20 and 0.80. To ensure stable mapping of various geometric shapes to a common generic domain and efficient learning of operators, the architecture of the proposed DNO was a combination of harmonic mapping and Fourier Neural Operator (FNO), where the generic domain was a square domain with a size of $\Omega_I = [0,1] \times [0,1]$ and a resolution of $128 \times 128$ (Fig. 2a). The neural operator was trained with 1000 samples generated using L2 loss function (more details provided in the Methods section).

The trained DNO exhibited stability in both prediction accuracy and generalization across various domains even with shapes and sizes different from the training dataset. Firstly, on test datasets ranging from 128 to 1024 resolution, the prediction error L2 remained stable (Fig. 2d), indicating that the diffeomorphism neural operator possessed mesh invariance. Furthermore, the model demonstrated remarkable predictive accuracy on larger domain sizes. Even when the domain size was scaled up to $30 \times 30$, the model could still accurately predict the solution function, with the prediction error negligible compared with the training dataset domain (Fig. 2b and Fig. 2e). When the model generalized (was applied) to domains with significant shape variations, such as hexagons and octagons, it shown robustness with domain changes, with prediction error L2 decreased by no more than 0.03 (Fig. 2c and Fig. 2f).

For the generalization performance of DNO on complex domains, we proposed an evaluation method called Domain Diffeomorphism Similarity (DDS). The relationship between DDS and prediction error on hexagon and octagon shape dataset is shown in Fig. 2g, which shows a strong correlation between the DDS and prediction error. When DDS is greater than 0.97, the L2 prediction error remains within 0.07, indicating stable prediction performance. Conversely, if DDS falls below this threshold, the prediction stability may decrease (Fig. 2g). It indicates that the DDS can be used as a mean to evaluate the generalization of diffeomorphism neural operators.

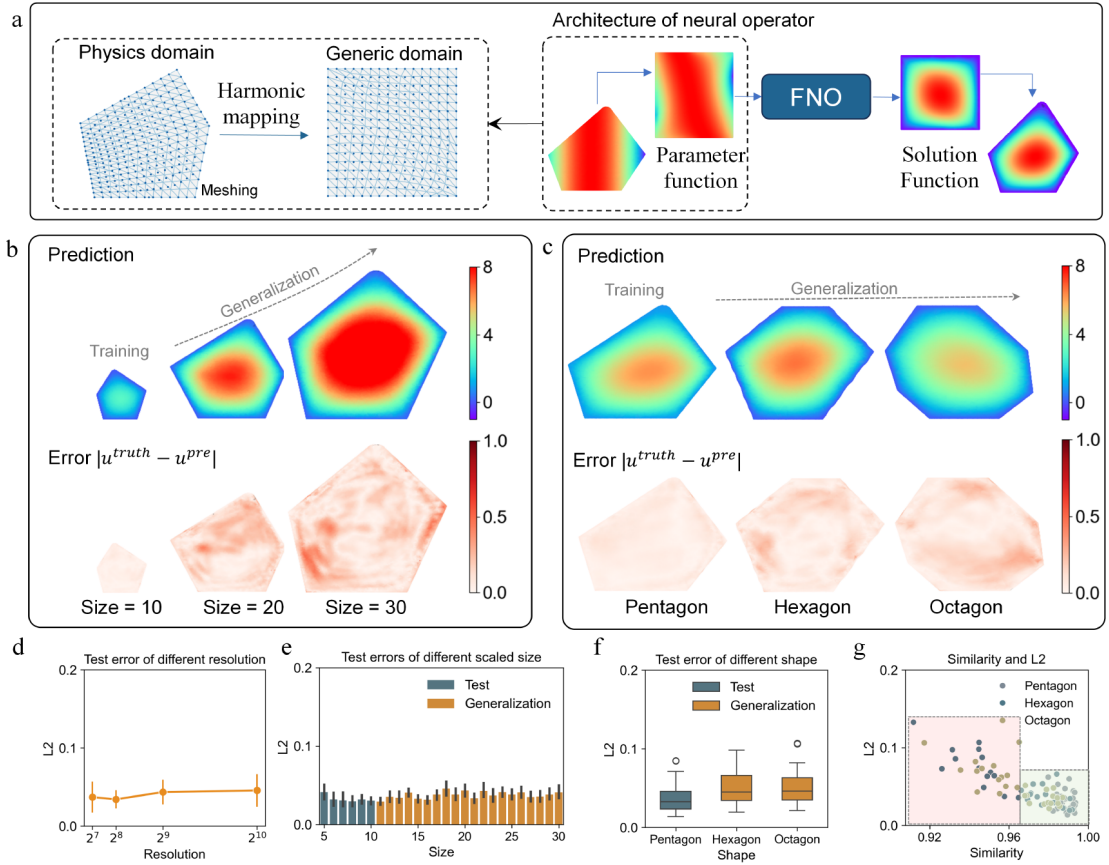

Fig. 2 | Application of diffeomorphism neural operator of Darcy flow cases on 2D domains. (a) The architecture of neural operator and the data process of diffeomorphic mapping and function sampling of parameter and solution function. (b) Generalization cases on pentagon domains with increased sizes. (c) An indication of generalization cases on domains with different shapes, including hexagon and octagon. (d) Test errors of different resolutions. (e) Test errors of different domain sizes ranging from 5 to 30. (f) Test errors of different shapes. (g) Relationship between geometric similarity on generic domain and prediction error.

**Fluid Dynamics equation**

Fluid dynamics analysis is important in various engineering disciplines, facilitating engineers and scientists in the design and optimization of a diverse array of systems in aerospace, automotive, energy and environmental infrastructures[32]. Most fluids can be modeled with the Navier-Stokes (N-S) equations[33] which are non-linear equations that describe the interplay of a velocity field $\boldsymbol{v}$ and a pressure field $p$ within a fluid domain $\Omega_i$:

$$\rho \left( \frac{\partial \boldsymbol{v}}{\partial t} + (\boldsymbol{v} \cdot \nabla) \boldsymbol{v} \right) = -\nabla p + \mu \Delta \boldsymbol{v} + \boldsymbol{f} \quad (5)$$

$$\nabla \boldsymbol{v} = 0 \quad (6)$$

where $\mu$ is the viscosity, $\rho$ is the fluid density and $\boldsymbol{f}$ is the body force. Here, the first $N$ velocity frames of the domain are deemed as the 'parameter functions', and the following velocity domain is the 'solution function'.

Traditional methods for solving fluid equations demand substantial computational resources and time. Each time the domain or parameters change, simulation must be rerun, which is a considerable obstacle for optimization challenges in engineering. To address this issue, our proposed DNO can swiftly generate solutions for flow fields across various domains and parameters. Our method was tested on two typical fluid problems: pipe flows (with different numbers of baffles) and airfoil flows (with different airfoils) on 2D domains.

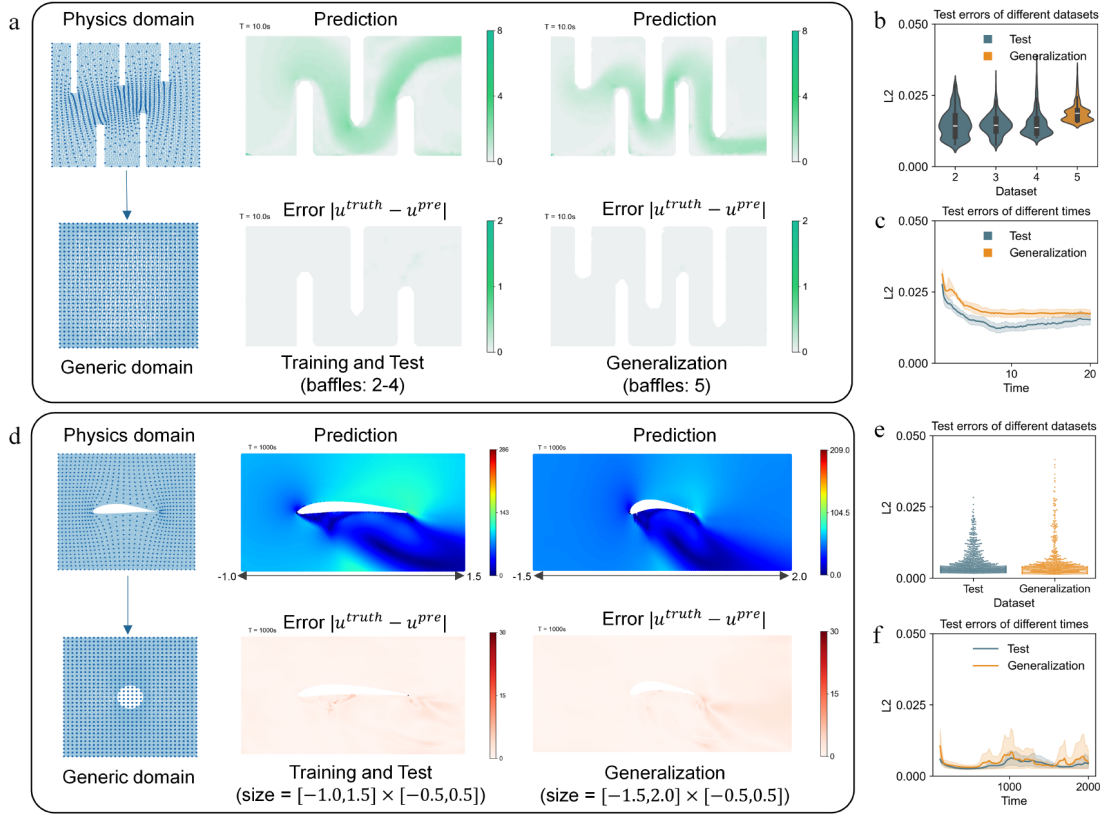

Fig. 3 | Application of diffeomorphism neural operator of flow dynamics cases on 2D domains. (a) Diffeomorphic mapping and test results of pipe flow cases. (b) Test errors of different dataset with different domains of pipe flow cases. (c) Temporal trend of error on test and generalization cases of pipe flow cases. (d) Diffeomorphic mapping and test results of airfoil flow cases. (e) Test errors of different datasets with different shapes and domains. (f) Temporal trend of errors on different datasets of airfoil flow cases.

For the pipe flow problem, the domain of the N-S equation consisted of the domains with varying numbers of baffles (range from 2 to 5) and changes in baffle positions and sizes. The initial velocity in each domain followed a random process, with values randomly selected between 1.0 and 4.0. Initially, different pipe geometries were mapped to a generic square domain (Fig. 3a) with a resolution of 128 ×128. Subsequently, the neural operator was trained on a dataset of 200 samples, where the number of baffles ranged from 2 to 4. It was then tested and generalized on a dataset where the number of baffles ranged from 2 to 5. Each category of baffles had 5 samples, and the generalization set each consisted of 5 domains with 5 baffles.

For the airfoil flow problem, the PDE equation was described by the Reynolds-Averaged Navier-Stokes (RANS) equations, which were derived from the more fundamental N-S equations by introducing additional terms to model turbulent effects. In this case, the varying domain with a size $[-1.0,1.5] \times [-0.5,0.5]$ comprised different airfoil shapes, while the parameters undergoing variation included the Reynolds number (range from 0.5 to 5 million) and the angle of attack (ranging of ±22.5 degrees). The airfoil shapes were obtained from the UIUC airfoil database[34,35]. Here, the domains containing airfoil shapes were diffeomorphically mapped into a square domain with a circular shape in the center (Fig. 3d), and the resolution was 128 ×128. The training dataset was generated randomly from the flow conditions described above, comprising a total of 100 samples of different airfoil shapes in the domain of size $[-1.0,1.5] \times [-0.5,0.5]$, while the test dataset was 20 different airfoil shapes and the generalization dataset was 10 airfoil shapes in a larger domain sized $[-1.5,2.0] \times [-0.5,0.5]$.

As mentioned in the Darcy flow cases, the architecture of DNO for pipe flow and airfoil flow was also the combination of harmonic mapping approach and FNO. During the training process, frames of the velocity fields from the preceding 5-time steps were utilized to predict the velocity field at the next time step. In the experiments involving both pipe flow and airfoil flow, our proposed framework demonstrated its ability to accurately predict solutions to complex fluid dynamics equations. Following training on a number of domains, it reliably predicted fluid behaviors across various domains and field parameters.

In the pipe flow experiments, prediction errors consistently remained under 0.04 for domains with 2 to 4 baffles in test datasets. Moreover, when DNO was generalized to domains with 5 baffles, the errors remained stable and also under 0.04 (Fig. 3b). It can be observed that the temporal errors in Fig. 3c, remained stability for both the test and generalization datasets, indicating the consistent performance across different times in various domains.

In the case of airfoils, similar results were observed. Prediction errors consistently remained below 0.05 across various airfoil shapes for both the test and generalization datasets (Fig. 3e). DNO also demonstrated consistent prediction accuracy over time, as depicted in both Fig. 3c and Fig. 3f. These findings suggested a promising approach for expediting the design process of pipe and airfoil shapes, particularly crucial for applications operating in complex environments that necessitate extensive testing and optimization, such as wide-speed domains. It is important to note that extensive testing across diverse parameters and shapes is still necessary.

**Mechanics equation**

The mechanical equation governing stress and deformation fields across various geometric domains is used in complex three dimensional (3D) mechanical scenarios in a wide range of engineering contexts, notably in component manufacturing, where material removal induces geometric changes, consequently impacting deformation fields amidst initial residual stress conditions. Here, the prediction of machining deformation of aircraft structural component was taken as an example, which was a problem of solving partial differential equations with variations in domain and parameter functions due to the material removal machining process. The geometric change caused by machining has a great influence on component deformation, rendering them as crucial as changes in the residual stress field for accurate deformation prediction. The deformation field of a component is affected by its internal residual stress field and the geometry of the component, and the relationship is governed by the following PDE, for $\Omega_i \in \mathbf{\Omega}$,

$$\nabla \sigma(x,y,z) + \nabla \left( c \nabla u_{u,\Omega_i}(x,y,z) \right) = 0 \quad (x,y,z) \in \Omega_i \tag{7}$$

$$u_{z,\Omega_i}(x,y) = 0 \qquad (x,y) \in E \tag{8}$$

where $\sigma(x,y,z)$ is the initial residua stress, which is the parameter function in this example. $\Omega_i \in \mathbf{\Omega}$ is the geometric domain of the component, which changes after each machining process. $E$ is the boundary condition points of the domain, which represents the constraint point of the component. $u_{z,\Omega_i}(x,y,z)$ is the deformation (solution) function of the component.

In this study, some typical aircraft structural component with pocket characteristics were chosen for testing, where the domains $\mathbf{\Omega}$ were a set of components with different geometries generated by machining processes. Because the geometric domains were 3D geometries with complex boundaries, the architecture of DNO was a combination of a volume parameterization (refer to Method section and Supplement III) of diffeomorphic mapping and FNO. The generic domain was a cuboid of size $\Omega_I = [0,600] \times [0,240] \times [0,30]$, as shown in Fig. 4b.

The data generation process included a random residual stress fields $\sigma(x,y,z) = c \times \sigma_0(x,y,z)$, where $c$ ranged from 0.85 to 1.15 and the function $\sigma_0(x,y,z)$ was the initial residual stress field of the component obtained by die forging simulation. There were two components with different geometries: Component 1 had 6 pockets, and Component 2 had 5 pockets, where the pockets were machined by material removal process along the Z axis, accomplished through 12 machining processes. Therefore, for each component, there were 12 domains generated by machining processes. The data generated on the first 10 domains of Component 1 were used as training dataset. The mean prediction error of the test dataset on the first 10 domains of Component 1 with randomly generated parameter function was 0.041, as shown in Fig. 4a and Fig. 4d. To assess the ability of generalization, the DNO predicted deformation solutions using the dataset of the last two domains of the 12 domains of Component 1 generated in the machining processes, yielding an L2 error of 0.044 (Fig. 4d).

Furthermore, the trained DNO on Component 1 was also applied to Component 2 (Fig. 4c). Because the deformation of components with distinctively different geometric structures exhibit significant variations, thus can have high challenge on the predictive capability of the neural operator. Through this example, we observed that the neural operator performed remarkably well in predicting deformation fields under substantial geometric variations, with an L2 prediction error of 0.051, as shown in Fig. 4d. The trend error of machining process was shown in Fig. 4e, which was stable in test and generalization dataset. Considering that even just the first 10 domains of the component were used for training, these results were acceptable for proving the quality of the proposed DNO.

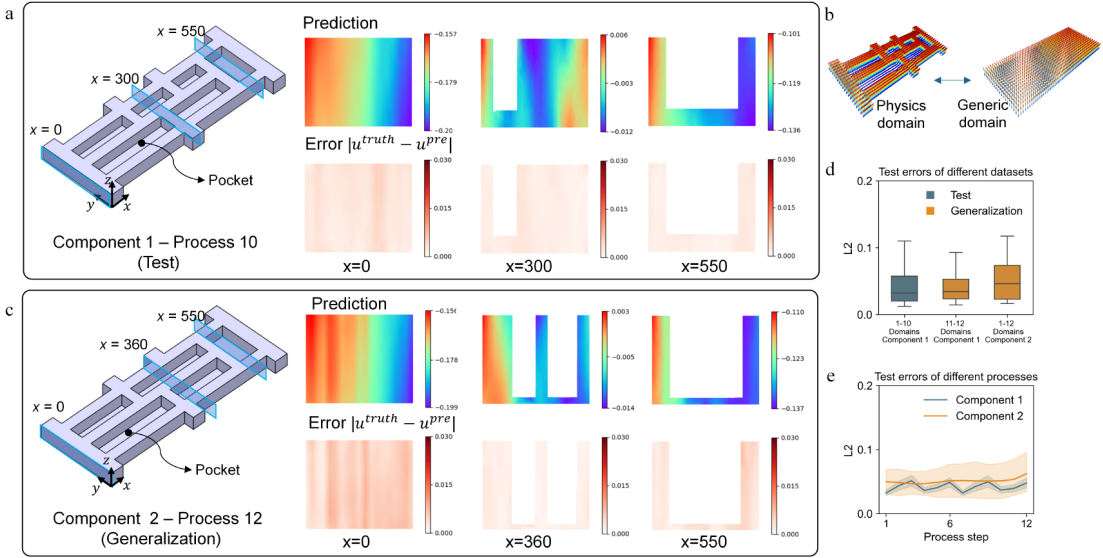

Fig. 4 | Application of DNO for deformation prediction cases on 3D domains. (a) Some test cases of prediction on Component 1 and Component 2. (b) Diffeomorphic mapping process from component domain to generic domain. (c) Test and generalization results on different datasets. (d) Test results along with the machining processes.

## Methods

### Diffeomorphism of various domains based on harmonic mapping

Here, the harmonic mapping used in 2D domain cases for establishing diffeomorphic mapping from physics domain $\Omega_i \in \boldsymbol{\Omega}$ to the generic domain $\Omega_I$ is summarized. Harmonic mapping is a useful tool for constructing diffeomorphisms under the conformal and bijective conditions (refer to reference[36] for additional details).

Given a physics domain $\Omega_i$ and the proposed generic domain $\Omega_I$, the mapping $\phi_i$ from $\Omega_i$ to $\Omega_I$ is harmonic if it satisfies equation $\Delta \phi_i = 0$. It is assumed that physics domain $\Omega_i$ has a boundary $B_p$, the generic domain $\Omega_I$ has a boundary $B_I$. Firstly, domain $\Omega_i$ is parameterized by meshing method (Fig. 5a), with vertexes $\boldsymbol{X}^{im} = (\boldsymbol{x}_1^{im}, \ldots, \boldsymbol{x}_n^{im})$, where $\boldsymbol{x}_i^{im}$ is the coordinates of the vector of mesh vertex. The aim of mapping $\phi_i$ is to find the corresponding coordinates of vertexes in the generic domain $\boldsymbol{X}^{Im}$. For the discrete meshing domain $\Omega_i$, the Laplacian operator $\Delta$ can be defined in the discrete form and represented by a matrix, and the equation can be rewritten as,

$$\boldsymbol{L}\boldsymbol{X}^{Im} = \boldsymbol{b} \tag{9}$$

$$L_{i,j} = \begin{cases} 1 & i = j \text{ and } i \in B_p \\ -w_{ij} & (i,j) \notin B_p \\ 0 & otherwise \end{cases} \tag{10}$$

where $\boldsymbol{L}$ is the Laplacian matrix of the domain meshes in $\Omega_i$, $\boldsymbol{X}^{Im}$ is an unknown vector of the of the mapping function at the vertexes of the mesh in $\Omega_I$, i.e., $\boldsymbol{X}^{Im} = (\boldsymbol{x}_1^{Im}, \ldots, \boldsymbol{x}_n^{Im})$, and $\boldsymbol{b}$ is the vector of vertexes on boundary $B_I$. $w_{ij}$ is the weight of the edge $e_{ij}$ connecting vertexes $i$ and $j$ in $\Omega_i$ and can be expressed by,

$$w_{ij} = \frac{1}{2} * (cotan\, \alpha_{ij} + cotan\, \beta_{ij}) \tag{11}$$

where $\alpha_{ij}$ and $\beta_{ij}$ are the two angles of edge $e_{ij}$.

To map the manifold $\Omega_i$ on to the same input manifold $\Omega_I$, the fixed boundary conditions (Dirichlet conditions) are used to solve the mapping solution, which imposes that the vectors within the domain are harmonic. Thus, the boundary vertexes of the physics domain are mapped to the boundary vertexes of the genetic domain. Then, considering the boundary conditions (vertexes of the boundary in the genetic domain) into $\boldsymbol{L}$ and $\boldsymbol{b}$, the modified equation is $\boldsymbol{L}_1\boldsymbol{X}^{Im} = \boldsymbol{b}_1$. After solving the linear functions, vectors $\boldsymbol{X}^{Im}$ of meshes in domain $\Omega_I$ could be obtained, as shown in Fig. 5a.

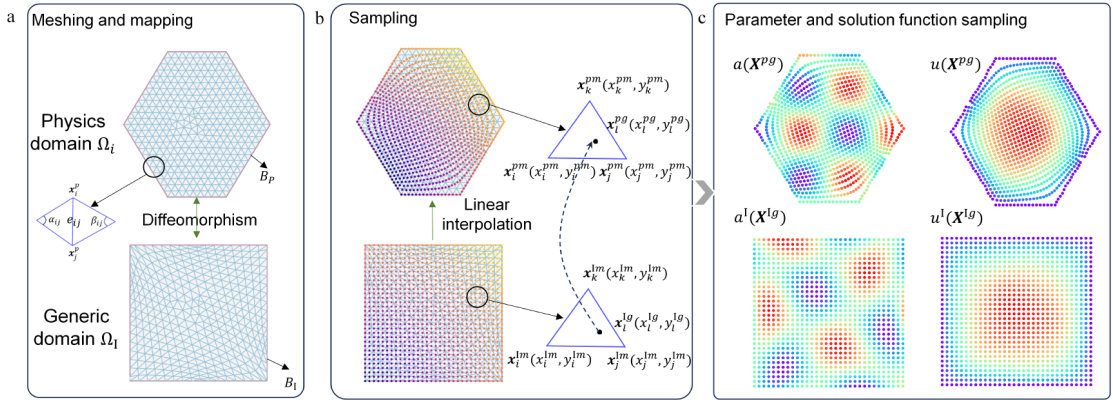

Fig. 5 | Diffeomorphism between physics domains and the generic domain based on harmonic mapping. (a) Calculation flow diagram of harmonic mapping. (b) Sampling method based on linear interpolation. (c) Sampling process of parameter and solution functions.

The parameter function and solution function are always represented by the values of sampled vectors. Harmonic mapping facilitates the mapping of various physics domains $\Omega_i$ onto the generic domain $\Omega_I$, each with distinct and irregular meshes. Therefore, to ensure the regular and uniform sampling of parameter function in the generic domain, the interpolation mapping method is used to get corresponding vectors and values in the physics domains.

For any sampling point in the generic domain obtained through uniform sampling, its coordinates can be linearly interpolated using the known vectors $X^{Im}$. Each vector $x_i^{Im} \in X^{Im}$ in the generic domain with coordinates. Therefore, to get the coordinates of the sampling points in the physics domain, the interpolation weights of the generic domain can be utilized. Taking the 2D plane as an example, suppose that the uniform grid point vector of the generic domain is $x_l^{Ig}(x_l^{Ig}, y_l^{Ig})$, the sampled vector of the physics domain after interpolation is $x_l^{pg}(x_l^{pg}, y_l^{pg})$. The calculation process of $x_l^{pg}(x_l^{pg}, y_l^{pg})$ is shown below,

$$\alpha = \frac{-(x_l^{Ig} - x_j^{Im})(y_k^{Im} - y_j^{Im}) + (y_l^{Ig} - y_j^{Im})(x_k^{Im} - x_j^{Im})}{-(x_i^{Im} - x_j^{Im})(y_k^{Im} - y_j^{Im}) + (y_i^{Im} - y_j^{Im})(x_k^{Im} - x_j^{Im})} \quad (12)$$

$$\beta = \frac{-(x_l^{Ig} - x_k^{Im})(y_i^{Im} - y_k^{Im}) + (y_l^{Ig} - y_k^{Im})(x_i^{Im} - x_k^{Im})}{-(x_j^{Im} - x_k^{Im})(y_i^{Im} - y_k^{Im}) + (y_j^{Im} - y_k^{Im})(x_i^{Im} - x_k^{Im})} \quad (13)$$

$$\gamma = 1 - \alpha - \beta \quad (14)$$

$$x_l^{pg} = \alpha \times x_i^{pm} + \beta \times x_j^{pm} + \gamma \times x_k^{pm} \quad (15)$$

$$y_l^{pg} = \alpha \times y_i^{pm} + \beta \times y_j^{pm} + \gamma \times y_k^{pm} \quad (16)$$

where $\alpha$, $\beta$ and $\gamma$ are interpolation weights calculated based on liner barycentric interpolation. $(x_i^{Im}, y_i^{Im})$, $(x_j^{Im}, y_j^{Im})$ and $(x_k^{Im}, y_k^{Im})$ are the mesh vertex vector of three points covering sampled grid points $x_l^{Ig}$. $(x_i^{pm}, y_i^{pm})$, $(x_j^{pm}, y_j^{pm})$ and $(x_k^{pm}, y_k^{pm})$ are the coordinates of the corresponding points of mesh triangle in the physics domain, as shown Fig. 5b.

After calculating all the sampling grid points $X^{pg}$ in the physics domain according to $X^{Ig}$ in the generic domain, the parameter function $a^I$ is mapped from $a^i$ in physics domain $\Omega_i$, and solution function $u^I$ is mapped from $u^i$ in physics domain $\Omega_i$, which can be expressed as $a^I(X^{Ig}) = a^i \circ \phi_i(X^{pg})$ and $u^I(X^{Ig}) = u^i \circ \phi_i(X^{pg})$ respectively in the generic domain, as shown in Fig. 5c. In neural operator learning process, the function $a^I(X^{Ig})$, together with $X^{pg}$ and $X^{Ig}$ are the input of the neural operator. The solution function $u^I(X^{Ig})$ is the output function of the neural operator. The sampling process is illustrated in Fig. 5b and Fig. 5c. It should be noted that the generic domain is the same as the output domain in this paper. In summary, through harmonic mapping and coordinate mapping, the non-uniformly sampled data in different physics domains can be uniformly used as the input of the subsequent neural operators in the generic domain.

**Diffeomorphism based on volume parameterization**

This study utilizes volume parameterization to establish the diffeomorphism from the physics domain of a 3D component to a generic cuboid domain. It should be noted that the generic domain can be other shapes because the DNO is flexible in using different shapes as the generic domain. The mapping process mainly includes the following steps: (1) geometric analysis of the component, (2) determination of the generic domain, and (3) volume parameterization. Finally, we provide a proof of the diffeomorphism of volume parameterization.

**Step (1):** Geometric analysis of the component. A die-forged structural component is selected in this study which has pocket features with the following characteristics: the geometric variations of the component only occur on the top surface of the component, and the cross-sections of the component along the length direction have a genus of 0.

**Step (2):** Determination of the generic domain. According to the geometric characteristics of the component, a cuboid domain is selected as the generic domain. In order to minimize the shape distortion of the component after mapping, a cuboid domain of 600mm × 240mm × 30mm is chosen as the generic domain.

**Step (3):** Volume parameterization. The top surface, bottom surface, and two side surfaces of the structural component are expressed through functions $T(x,y)$, $B(x,y)$, $L(x,z)$, and $R(x,z)$, respectively. Then volume parameterization can be expressed by the following formula,

$$\begin{bmatrix} x_l^{lg} \\ y_l^{lg} \\ z_l^{lg} \end{bmatrix} = f\left(\begin{bmatrix} x_l^{pg} \\ y_l^{pg} \\ z_l^{pg} \end{bmatrix}\right) = \begin{bmatrix} x_l^{pg} \\ \frac{240 \times \left(y_l^{pg} - R(x_l^{pg})\right)}{L(x_l^{pg}) - R(x_l^{pS})} \\ \frac{30 \times \left(z_l^{pS} - B(x_l^{pS}, y_l^{pg})\right)}{T(x_l^{pg}, y_l^{pg}) - B(x_l^{pg}, y_l^{pg})} \end{bmatrix}; \begin{bmatrix} x_l^{lg} \\ y_l^{lg} \\ z_l^{lg} \end{bmatrix} \in \Omega_I, \begin{bmatrix} x_l^{pg} \\ y_l^{pg} \\ z_l^{pg} \end{bmatrix} \in \Omega_i \quad (17)$$

where $f$ represents the volume parameterization function, which is the mapping from the physics domain of the component to the generic cuboid domain, $(x_l^{lg}, y_l^{lg}, z_l^{lg})$ denotes the coordinates of a point in the generic domain after mapping, and $(x_l^{pg}, y_l^{pg}, z_l^{pg})$ represents the coordinates of a point in the physics domain of the component before mapping. The method can achieve mappings from domains of similar components to the generic domain. The proof of the diffeomorphism of volume parameterization is provided in Supplement. III.

**Diffeomorphism neural operators using Fourier integral**

In the proposed DNO framework, after obtaining an equal number of sampling points in the generic domain, to ensure the match between the dimensions of the parameter functions and solution functions from different domains, the selection of the framework of the neural operator can be flexible, such as DeepONet or FNO. In this study, we used FNO and incorporated geometric information from both the physics domain and the generic domain into the Fourier layer. Here, we provide a brief overview of our architecture.

FNO consists of the N-layer nonlinear operator layer $K$, a lifting operator $P$ and a local projection operator $Q$, which is defined as follows,

$$\mathcal{G}_\theta := Q \circ \mathcal{L}_N \circ \cdots \mathcal{L}_n \cdots \circ \mathcal{L}_1 \circ P \quad (18)$$

where $\mathcal{G}_\theta$ is the neural operator, $P$ is the lifting operator that encodes the lower dimension input function $(a^I, X^{pg}, X^{lg})$ in the generic domain into higher dimensional space represented by $z_0$, and $Q$ is the projection operator that maps the last representation $z_N$ on to the output domain. Generally, both the lifting operator $P$ and the projection operator $Q$ are fully neural networks. The N-layer nonlinear operator layer $\mathcal{L}_1, \cdots, \mathcal{L}_N$ has the same structure as follows,

$$\mathcal{L}_n(z_n) = \sigma(W_n(z_n) + b_n + K_n(z_n)) \tag{19}$$

where $\sigma: \mathbb{R} \to \mathbb{R}$ is a point-wise non-linear activation function, the weight matrix $W_n$ and the bias $b_n$ define an affine pointwise map $W_n(z_n) + b_n$. $K_n(z_n)$ which are integral kernel operators defined as,

$$K_n(z_n) = \mathcal{F}^{-1}(R_n \cdot \mathcal{F}(z_n)) \tag{20}$$

where $R_n$ is the Fourier-domain weight matrix, $\mathcal{F}$ and $\mathcal{F}^{-1}$ are Fast Fourier transform and inverse Fast Fourier transform respectively.

In the diffeomorphism neural operators, we have improved the structure of Fourier neural operators by incorporating geometric information from both the physical domain and the generic domain within the Fourier layer. We have transformed the Fourier operator layer into:

$$\mathcal{L}_n(z_n, X^{pg}, X^{lg}) = \sigma\left(W_n(z_n) + G^p(X^{pg}) + G^I(X^{lg}) + b_n + K_n(z_n)\right) \tag{21}$$

where $G^p$ and $G^I$ are weight matrix used to map geometry information from the physics and the generic domain, respectively. Together with $W_n$ and $b_n$, they constitute an affine pointwise map. Then the components of diffeomorphism of the physics domain and the generic domain are combined with the FNO component, forming the diffeomorphism neural operator (DNO), which could be adapted to different physics domains.

The neural operator $\mathcal{G}_\theta$ could be defined as to achieving the map between the input function to the output function, i.e., $\mathcal{G}_\theta: (a^I, X^{pg}, X^{lg}) \to u^I$, where $\mathcal{G}_\theta$ is parametrized by $\theta$. The learning objective can be defined as to minimizing the empirical loss on dataset $\{a^I, X^{pg}, X^{lg}, u^I\}$. The loss function is based on L2 Relative error, as shown below,

$$Loss = E_{\Omega_i \sim \Omega} \frac{\|\mathcal{G}_I - \mathcal{G}_\theta\|^2}{\|\mathcal{G}_I\|^2} = \frac{1}{N*M} \sum_{i=1}^{N} \sum_{j=1}^{M} \frac{\|u_i^{truth} - u_i^{pre}\|^2}{\|u_i^{truth}\|^2} \tag{22}$$

where $u_i^{truth}$ and $u_i^{pre}$ are predicted and truth solution function sampled by $X^{pg}$ (total M sampled points) in the $i^{th}$ domain (total $N$ domains).

**Evaluating the generalization ability of diffeomorphism neural operators**

One of the primary goals of most machine learning algorithms is to generalize the algorithms to be applicable to new data which are not in the training dataset. For the proposed DNO, ensuring and quantifying generalization is a challenging issue. In machine learning, the evaluation of generalization ability concerns the similarity between data. Here, we focus on the generalization of neural operators across different domains. The proposed DNO learns mappings between parameter functions and solution functions in the generic domain. According to the sampling process mentioned in previous sections, both parameter functions $a(X^{pg})$ and solution functions $u(X^{pg})$ can be regarded as mappings on geometric parameters $X^{pg}$ on the physics domain. Therefore, the generalization ability of the DNO is mainly influenced by the geometry domains. Consequently, evaluating the similarity between different domains is a direct measure for assessing the generalization of the neural operator. Here, we introduce a measure called domain diffeomorphism similarity (DDS) to evaluate the generalization of diffeomorphism neural operators on new domains. Specifically, the DDS evaluates the similarity between geometries mapped from two different domains in a generic domain.

Considering the two diffeomorphisms $\phi_i$ and $\phi_j$, from domains $\Omega_i$ and $\Omega_j$ to the generic domain $\Omega_I$, the domain parametrization functions are $\boldsymbol{X}_i^{pg}$ and $\boldsymbol{X}_j^{pg}$. The aim is to evaluate the similarities between $\boldsymbol{X}_i^{pg}$ and $\boldsymbol{X}_j^{pg}$. Actually, $\boldsymbol{X}_i^{pg}$ and $\boldsymbol{X}_j^{pg}$ could be regarded as two geometric images [36] in the generic domain space. Therefore, the Normalized Cross-Correlation (NCC) can be used to evaluate the similarity between the two geometric images with the size of $M \times N \times C$ (For a 2D domain with a resolution of 128, the size is $128 \times 128 \times 2$), which can be expressed as,

$$NCC\left(\boldsymbol{X}_i^{pg}, \boldsymbol{X}_j^{pg}\right) = \frac{\sum_{l=1}^{MN}\sum_{c=1}^{C}\left(x_{i,n,c}^{pg} - \bar{x}_{i,c}^{pg}\right)\left(x_{j,n,c}^{pg} - \bar{x}_{j,c}^{pg}\right)}{\sqrt{\sum_{l=1}^{MN}\sum_{c=1}^{C}\left(x_{i,n,c}^{pg} - \bar{x}_{i,c}^{pg}\right)^2}\sqrt{\sum_{l=1}^{MN}\sum_{c=1}^{C}\left(x_{j,n,c}^{pg} - \bar{x}_{j,c}^{pg}\right)^2}} \tag{23}$$

where $\bar{x}_{i,c}^{pg}$ and $\bar{x}_{j,c}^{pg}$ are the average values of tensors $\boldsymbol{X}_i^{pg}, \boldsymbol{X}_j^{pg}$ respectively on dimension $c$. This formula describes the similarity between two domains with multiple dimensions. In order to evaluate whether the trained operator is generalized to be applicable to new domains that were not in the training dataset, the NCC between the new domain and the domain in the training set is calculated and the average value is taken as the measure, i.e., DDS, to evaluate the generalization ability of DNO for this new domain. The results of DDS and generalization error are shown in Fig. 2f. It can be seen that there is a strong correlation between DDS and prediction error, indicating that DDS can be used as a measure to evaluate the generalization ability of diffeomorphism neural operators.

**Conclusions and further work**

Although neural operators based on deep learning have promising potential to significantly increase the computation efficiency of traditional PDE solving, the neural operators reported by previous researchers only solved PDEs on fixed domains - this was a major limitation in many scientific and engineering applications. Therefore, the main contribution of this research is a proposed general neural operator learning method which can solve PDEs on various domains. In principle, the proposed diffeomorphism neural operator (DNO) learns the physics law governed by a PDE on a generic domain (a new concept proposed in this research) which is a diffeomorphism with all physics domains governed by the same PDE. The learning ability of the proposed DNO had been proved and validated through experiments on both static and dynamic PDEs, including Darcy flow equation, N-S equation and mechanics equation. The generalization ability of the proposed NDO on new domains that are different from the training dataset had also been evaluated, including variations in both size and shape. The experimental results proved that the method exhibited strong operator learning ability, particularly on complex various domains.

In the literature review, we have not seen methods like the proposed DNO framework which could solve PDEs on various domains, thus DNO has the potential to be used in a wider range of applications in science and engineering, such as optimization in design, manufacturing process and scientific experiments. Our future research directions include optimization of the DNO framework by enhancing the learning ability of operators through introducing new methods based on spatial invariance, rotational invariance, and other neural network techniques.

**Data and Code Availability**

The code and dataset scripts used in this study are available on the GitHub repository https://github.com/Zhaozhiwhy/Diffeomorphism-Neural-Operator.

**Acknowledgements**

This research was funded by the National Science Fund of China for Distinguished Young Scholars under Grant 51925505 (Yingguang Li) and National Natural Science Foundation of China under Grant 52175467 (Changqing Liu).

**Author Contributions Statement**

Zhiwei Zhao: Methodology, Conceptualization, Formal analysis, Writing - Original draft & editing.
Changqing Liu: Methodology, Conceptualization, Writing - review & editing, Funding acquisition.
Yingguang Li: Methodology, Conceptualization, Funding acquisition, Supervision.
Zhibin Chen: Formal analysis, Writing - Original draft.
Xu Liu: Methodology, Conceptualization, Writing - review & editing.

**Competing interests**

The authors declare no competing interest.